\documentclass[12pt,reqno]{amsart}
\usepackage{verbatim}
\usepackage{amscd}
\usepackage{amsfonts}
\usepackage{epsfig}

\begin{document}

% The following describe useful names and variations.
\newcommand{\Bc}{{\mathcal B}}
\newcommand{\Ec}{{\mathcal E}}
\newcommand{\Hc}{{\mathcal H}}
\newcommand{\Kc}{{\mathcal K}}

\newcommand{\Lc}{{\mathcal L}}

\newcommand{\Om}{\Omega}
\newcommand{\pal}{\partial}

\newcommand{\Rb}{\mbox{$\mathbb R$}}
\newcommand{\Sc}{{\mathcal S}}
\newcommand{\Sct}{\tilde{\Sc}}

\newcommand{\Sg}{\Sigma}

\newcommand{\ubf}{\bf{u}}
\newcommand{\vap}{\varphi}
\newcommand{\vbf}{\bf{v}}

\newcommand{\wbf}{\bf{w}}

\newcommand{\zbf}{\bf{z}}
\newcommand{\Wc}{{\mathcal W}}

% The following establish some abbreviations for mathematical symbols.
% 
\newcommand{\bdy}{\partial \Omega}
\newcommand{\bGo}{\partial G_1}
\newcommand{\bGl}{\partial \Gal}

\newcommand{\delj}{\delta_j}
\newcommand{\delk}{\delta_k}
\newcommand{\delm}{\delta_m}
\newcommand{\delo}{\delta_1}

\newcommand{\etab}{\overline{\eta}}
\newcommand{\etaj}{\eta_j}

\newcommand{\Gal}{G_{\alpha}}

\newcommand{\hav}{\overline{h}}
\newcommand{\hjhat}{\hat{h_j}}

\newcommand{\nualp}{\nu \alpha}
\newcommand{\nuj}{\nu_j}
\newcommand{\nuk}{\nu_k}
\newcommand{\Omb}{\overline{\Omega}}

\newcommand{\psit}{\tilde{\psi}}

\newcommand{\RN}{{\Rb}^N}
\newcommand{\Rtwo}{{\Rb}^2}
\newcommand{\Rte}{{\Rb}^3}

\newcommand{\Scone}{\Sc_1}

\newcommand{\Sgb}{\overline{\Sg}}
\newcommand{\Sgz}{\Sigma_0}
\newcommand{\sinc}{\mbox{sinc}}
\newcommand{\sj}{s_j}
\newcommand{\sjt}{\tilde{\sj}}
\newcommand{\sk}{s_k}
\newcommand{\skt}{\tilde{\skt}}
\newcommand{\soj}{\st_{1j}}

\newcommand{\stj}{\st_{2j}}
\newcommand{\st}{\tilde{s}}
\newcommand{\tsoj}{\widetilde{s_{1j}}}
\newcommand{\tstj}{\widetilde{s_{2j}}}

\newcommand{\vaph}{\hat{\varphi}}
\newcommand{\vapt}{\tilde{\vap}}

\newcommand{\vbfm}{\vbf_m}
\newcommand{\vbh}{\hat{\vbf}}
\newcommand{\vbt}{\tilde{\vbf}}
\newcommand{\vnu}{\vbf \cdot \nu}
\newcommand{\vnum}{\vbf_m  \cdot \nu}

\newcommand{\wnu}{\wbf \cdot \nu}

\newcommand{\vh}{\hat{v}}

\newcommand{\wbj}{\wbf_j}
\newcommand{\wbk}{\wbf_k}
\newcommand{\wt}{\tilde{\wbf}}
\newcommand{\wbtj} {\wt_j}
\newcommand{\wbtk}{\wt_k}

% the following macros are for common words and combinations.

\newcommand{\foral}{\quad \mbox{for all} \quad}

\newcommand{\wrt}{\ \mbox{with respect to} \ }

\newcommand{\xand}{\quad \mbox{and} \quad }
\newcommand{\xfor}{ \quad \mbox{for} \ }
\newcommand{\xiff}{\ \mbox{if and only if} }
\newcommand{\xxiff}{\ \mbox{if and only if} \ }

\newcommand{\xhSfn}{\ \mbox{harmonic Steklov eigenfunction}}
\newcommand{\xSfn}{\ \mbox{Steklov eigenfunction} }
\newcommand{\xSvl}{\ \mbox{Steklov eigenvalue}}

\newcommand{\xon}{\qquad \mbox{on} \quad}
\newcommand{\xor}{\qquad \mbox{or} \quad}
\newcommand{\xwhen}{\qquad \mbox{when} \quad}
\newcommand{\xwith}{\qquad \mbox{with} \quad}

\newcommand{\xwlog}{\ \mbox{without loss of generality} \ }

% The following are macros for mathematical operations and combinations.
%
\newcommand{\ang}[1]{\langle#1\rangle}
\newcommand{\curl}{\mathop{\rm curl}\nolimits}

\newcommand{\deq}{:= }
\newcommand{\deqs}{\ :=\ }
\newcommand{\Dnu}{D_{\nu}}
\newcommand{\Div}{\mathop{\rm div}\nolimits}
\newcommand{\dsg}{\ d \sigma}
\newcommand{\dtx}{\ d^2 x}

\newcommand{\dv}{\ d^n x}

\newcommand{\eqs}{\ =\ }

\newcommand{\geqs}{\ \geq \ }
\newcommand{\gradchi}{\nabla \chi}
\newcommand{\gradphi}{\nabla \varphi}
\newcommand{\grads}{\nabla s}
\newcommand{\gradu}{\nabla u}
\newcommand{\gradv}{\nabla v}

\newcommand{\Iby}{\int_{\bdy} \ }
\newcommand{\IMby}{|\bdy|^{-1} \, \int_{\bdy} \ }

\newcommand{\IbG}{\int_{\pal \Gal} \ }
\newcommand{\IGal}{\int_{\Gal} \ }
\newcommand{\IOm}{\int_{\Om} \ }
\newcommand{\IOmL}{\int_{\Om_L} \ }

\newcommand{\leqs}{\ \leq \ }
\newcommand{\mns}{\ - \ }
\newcommand{\nm}[1]{\left\Vert#1\right\Vert}
\newcommand{\Omo}{{\Om}_1 }
\newcommand{\OmL}{{\Om}_L }

\newcommand{\pls}{\, + \, }
\newcommand{\plms}{\ + \ }

% The following are TeX operational abbreviations.
\newcommand{\barr}{\begin{eqnarray}}
\newcommand{\beq}{\begin{equation}}
\newcommand{\bpf}{\begin{proof} \quad }
\newcommand{\btm}{\begin{thm} \quad }

\newcommand{\earr}{\end{eqnarray}}
\newcommand{\eeq}{\end{equation}}
\newcommand{\epf}{\end{proof}}
\newcommand{\etm}{\end{thm} \quad }

\newcommand{\ra}{\rightarrow}

%  The following are the theorem like structures used here.
%
\newtheorem{thm}{Theorem}[section]
\newtheorem{cor}[thm]{Corollary}
\newtheorem{cond}{Condition}
\newtheorem{lem}[thm]{Lemma}
\newtheorem{prop}[thm]{Proposition}
\numberwithin{equation}{section}
\renewcommand{\theequation}{\thesection.\arabic{equation}}
\setcounter{secnumdepth}{2}
\newcommand{\cnd}[2]{\par\medskip\noindent{\bf Condition #1}. {\em
 #2}\par\medskip}

%
% The following are various spaces of functions used
%
\newcommand{\Cc}{C_c^1 (\Om)}
\newcommand{\Cz}{C^0 (\Omb)}
\newcommand{\Ctest}{C_c^{\infty} (\Om)}

\newcommand{\Harm}{\Hc^1(\Om)}
\newcommand{\Harmo}{\Hc^1(\Omo)}
\newcommand{\Halfby}{H^{1/2}(\bdy)}
\newcommand{\Hone}{H^1(\Om)}
\newcommand{\HSgz}{H^1_{\Sg0}(\Om)}
\newcommand{\Hz}{H_0^1(\Om)}

\newcommand{\Lpby}{L^p (\bdy; \dsg)}
\newcommand{\Ltby}{L^2 (\bdy, \dsg)}
\newcommand{\LtSg}{L^2 (\Sigma, \dsg)}
\newcommand{\LmSg}{L^2_m (\Sigma, \dsg)}
\newcommand{\Lt}{L^2 (\Om)}

\newcommand{\Wp}{W^{1,p}(\Om)}

\title[Central Approximation of Harmonic Functions.] 
{Boundary Integrals and  Approximations of \\   Harmonic Functions}
\author[Auchmuty \& Cho]{Giles Auchmuty and Manki Cho}

\address{Department of Mathematics\\
    University of Houston\\
    Houston, Tx 77204-3008, USA. }
\email{auchmuty@uh.edu \xand realmann@math.uh.edu}

\date{\today}

\begin{abstract} 
Steklov expansions for a harmonic function on  a rectangle are analyzed.
The   value of a harmonic function at the center of a rectangle is shown to be well approximated by the mean value of 
the function on the boundary plus a very small number  (often 3 or fewer) of additional boundary  integrals.
Similar approximations  are found for the central values  of solutions of Robin and Neumann  boundary value problems.
These results are based on finding explicit expressions for the Steklov eigenvalues and eignfunctions. \end{abstract}
\maketitle

\section{Introduction}\label{s1}

One of the best known theorems about solutions of Laplace's equation is the mean value theorem that says the value of a harmonic function $h$ at the center of the ball is the mean 
value of the boundary values of $h$ on the boundary.
Epstein \cite{Ep} proved that if this holds for all harmonic functions on a region $\Om$ then $\Om$ must be a ball.

In a number of papers (see \cite{Au1} - \cite{Au4}), the first author has studied the representation of solutions of 
elliptic boundary value problems via orthogonal expansions involving  the Steklov eigenfunctions. 
These representations for harmonic problems may be regarded as   spectral representations of the Poisson integral  operator
for the solution of the harmonic Dirichlet problem on $\Om$.   

Steklov expansions give  the value of a harmonic function at any point $x \in \Om$ is given by a series whose first term is the mean value of the boundary values of the function.
See \eqref{e3.9} below. The coefficient in these expansions only involve the boudnary values of the functions and the Steklov eigenfunctions as described in \eqref{e3.11} since the 
definition of a Steklov eigenfunction relates an inner product over the domain to the $L^2$- inner product on the boundary. Theorem \ref{T3.3} below
that these Steklov expansions converge both in $H^1$ and uniformly on compact subsets of the region.

Here the special case of harmonic functions on a rectangle in the plane is investigated and some strong convergence properties are found. First the Steklov eigenvalues and eigenfunctions
for a rectangle are described explicitly in section 4. The eigenfunctions are classified by their symmetries and by the determining equation for a paprameter related to the eigenvalues.
The results hold for arbitrary rectangles in the plane upon scaling and possibly a rotation. In their recent survey, Girouard and Polterovich \cite{GP} outline similar results for the case of a square.
 
In the doctoral thesis of the second author \cite{Cho}  many aspects of the convergence of  Steklov approximation of harmonic functions on rectangles were studied numerically.
In particular Cho observed that these coefficients dacayed to zero very rapidly and that these Steklov approximations converged very rapidly away from the boundary of the rectangle.

The analysis here shows that this decay rate actually is exponential as described in theorems \ref{T5.1} and  \ref{T6.1}. The usual mean value theorem for discs and balls in 2 and higher dimensions
reflects the fact that all other Steklov eigenfunctions, (besides the first one which is constant,) are zero at the origen. Thus a natural question is what are the corresponding formulae at the center of
a rectangle or square?

At the origin 3 out of the 4 families of Steklov eigenfunctions are zero and the remaining family involves functions that are
exponentially  small at the center. 
As a consequence, the Steklov approximations converge geometrically to the true value of the harmonic function at the
center.
The usual mean value theorem for discs and balls in 2 and higher dimensions reflects the fact that all other  Steklov eigenfunctions, 
( besides the first one which is constant,)  are zero at the origin. 
Our observation is that, on a square, the use of just one  further boundary integral is sufficient to guarantee relative accuracy of
better than 4\% for the central value. 
Specific error estimates in terms of the 2-norm of the data is provide in Theorem \ref{T5.2} for a square and from a similar analysis 
involving theorem \ref{T6.1} for a rectangle.

In section 7, analogous formulae for the central values of solutions of Robin and Neumann boundary value problems are 
derived.
These involve the same  boundary integrals as for  the Dirichlet problem but with new weightings depending on the boundary
condition. 
In particular these formulae indicate how solutions of Robin problems with non-zero mean values for the boundary data
behave as it converges to a Neumann problem.

A referee has pointed out that some other exponential convergence results are described in  theorem 1.1 of Hislop and
 Lutzer \cite{HL}.
 The results here differ from those of \cite{HL} as this paper studies  pointwise estimates  at specific points - not $H^1$ bounds. 
Also here our  boundary  is Lipschitz but not $C^1$ so the very general methods used in \cite{HL} are not applicable.  
The recent review paper \cite{GP} discusses a number of features of Steklov eigenproblems that change significantly 
when the boundary is no longer smooth.

\vspace{1em}

%%% New Section 2 ****

\section{Notation and Definitions.}\label{s2}

Let $\Om$ be a non-empty, bounded, connected, open subset of $\RN, \ n \geq 2$, with boundary  $\bdy$.  
Such a set $\Om$ is called a {\it region} and let its boundary be $\bdy$.
Let $\Lt$ be the usual real Lebesgue space of all functions $\vap: \Om \ra [-\infty, \infty]$
which are square integrable with respect to Lebesgue measure on $\Om$.

For our analysis we also require some mild regularity conditions on $\Om$ and $\bdy$. 
Our basic criteria is that the boundary $\bdy$ is Lipschitz in the sense of Evans and Gariepy,
\cite{EG}, section 4.2.  

Let $\sigma, \dsg$ represent Hausdorff $(n-1)-$dimensional measure and integration  with respect to this measure  respectively.  
When $\bdy$ is Lipschitz, the unit outer normal $\nu(x)$ is defined  $\sigma \, a.e.$ on $\bdy$  and the  boundary traces of 
functions in $\Wp$ are well-defined. 
The trace operator will be noted $\gamma$ and usually will be omitted in boundary integrals. 
See section 4.3 loc. cit. In particular, the divergence  theorem holds in the form
\beq \label{e2.3}
\IOm D_j \, \psi \dv \eqs \Iby \psi \, {\nu}_j \dsg \qquad \mbox{for all} \ 1\leq j \leq n, \ \psi \in W^{1,1}(\Om).
\eeq
This follows from theorem 1, section 4.3 of \cite{EG} with p =1. 
Our requirement is that $\Om$  satisfy 

\noindent{\bf Condition A:}  \quad  {\it $\Om$ is a bounded open region and $\bdy$ is a finite union of Lipschitz surfaces
with finite surface area. }

When (A) holds and $\psi \in \Hone$, then the standard trace theorem implies that $\gamma(\psi) \in \Lpby$
for $p \in [1,p_T]$ where $p_T = 2(n-1)/(n-2)$ when $n \geq 3$ and for all $p \in [1,\infty)$ when $n=2$.
In particular they  are in $\Ltby$ for every n.

Let $ \Ltby$ be the usual space of   $L^2-$functions on $\Sigma, \bdy$ respectively.
The inner product and norm on $\Ltby$ is   
\beq \label{e2.5}
\ang{u,v}_{\bdy} \deqs \IMby u \, v \dsg.  \eeq

When $p=2$, the space $W^{1,2}(\Om)$ will be denoted  $\Hone$ and an equivalent inner product is
\begin{equation}\label{ip2}
{[u , v]}_{\pal}  \deq \IOm \  \gradu  \cdot \nabla v \  dx \pls
{|\bdy|}^{-1} \, \int_{\bdy} \  u \; v \ \dsg.
\end{equation}
The corresponding norm will be denoted by  $ {\nm{u}}_{\pal}$.
The proof that this  norm is equivalent to the usual $(1,2)-$norm on $\Hone$ when
(A) holds is Corollary 6.2 of \cite{Au1} and also is part of theorem 21A of \cite{Z}. 

A function $u \in \Hone$ is said to be {\it harmonic on $\Om$} provided it is a solution of 
Laplace's equation in the usual weak sense. Namely
\beq \label{Wharm}
\IOm \ \gradu \cdot \nabla \vap \ dx \eqs 0 \qquad \foral \ \vap \in \Cc.
\eeq
Here $\Cc$ is the set of all $C^1-$functions on $\Om$ with compact support in $\Om$.

Define $\Harm$ to be the space of all such $H^1-$harmonic functions on $\Om$. 
When (A)  holds, the closure of $\Cc$ in the $H^1-$norm is the usual Sobolev space $\Hz$.
Then (\ref{Wharm}) is equivalent to saying  that $\Harm$ is $\pal-$orthogonal to 
$\Hz$. This may be expressed as
 \beq \label{e2.7}
\Hone \eqs \Hz \oplus_{\pal}  \Harm,
\eeq
where $\oplus_{\pal}$ indicates that this is a $\pal-$orthogonal decomposition.
This result  is also discussed in section 22.4 of \cite{Z}.

%% need to define inner product on H^1, harmonic function on Omega. The subspace \Hc. 
The terminology of Axler, Bourdon and Ramey \cite{ABR} for harmonic functions will be used in discussing 
harmonic functions. Other standard notation from the theory of elliptic boundary value problems should be 
taken as that used  in Chipot \cite{Ch}.

\vspace{1em}
 %  %%%%%
% ####### Section 3 ######
% This is line 320 approx
\section{Steklov Approximations of Harmonic functions. } \label{s3}

Let $\Om$ be a bounded region in $\RN$ that satisfies (A).
It is well-known that the Bergman space of all real $L^2-$harmonic functions on $\Om$ is a reproducing
kernel Hilbert spaces (RKHS) with respect to the usual $L^2-$ inner product.
Since we wish to prove results involving boundary integrals, we will concentrate on the subspace
$\Harm$ of all $H^1-$ harmonic functions. 
This is also a RKHS of functions on $\Om$ but now the functions have boundary traces - see Auchmuty
\cite{Au2} and  \cite{Au3} for results about such  spaces. 
Trace theorems imply that for regions satisfying (A) the traces of functions in $\Harm$ will be in $\Ltby$,
so  the integrals involved here are  well-defined and finite.  

 A non-zero function $s \in \Hone$  is said to be a {\it harmonic Steklov eigenfunction}  on $\Om$ 
 corresponding to the  Steklov eigenvalue $\delta$ provided $s$ satisfies
\beq \label{HSeqn}
\IOm \;  \nabla s  \cdot  \nabla v  \; dx \eqs \delta  \   \ang{s,v}_{\bdy} \eqs \delta \IMby s \, v \dsg.  \qquad 
\mbox{for all} \quad v \in \Hone.
\eeq

  This is the weak form of the boundary value problem
\beq \label{Weeqn}
\Delta \, s \eqs 0 \qquad \mbox{on $\Om$ with} \quad  \Dnu \, s \eqs \delta \ |\bdy|^{-1} \  s 
\quad \mbox{on} \ \bdy.
\eeq
Here $\Delta$ is the Laplacian and $\Dnu \, s(x) := \nabla s(x) \cdot \nu(x)$ is the unit 
outward normal derivative of s    at a point on the boundary. 

Descriptions of the analysis of these eigenproblems may be found in Auchmuty \cite{Au1} - \cite{Au4}.
These eigenvalues and a corresponding family of $\pal-$orthonormal eigenfunctions 
may be found using variational principles as described in sections 6 and 7 of Auchmuty \cite{Au1}. 
$\delta_0 = 0$ is the least eigenvalue of this problem corresponding to the eigenfunction
$s_0 (x) \equiv 1$ on $\Om$. This eigenvalue is simple as $\Om$ is connected. 
Let the first k Steklov eigenvalues be $0 = \delta_0 < \delta_1 \le \delta_2  \leq \ldots \leq \delta_{k-1}$ 
and $s_0, s_1, \ldots, s_{k-1}$ be a corresponding  set of  $\pal-$orthonormal eigenfunctions. 
The k-th eigenfunction $s_k$ will be a maximizer of the functional
\beq \label{e3.3}
\Bc(u) \deqs \IMby  \, | \gamma u|^2 \, \dsg
\eeq
over the subset $B_k$ of functions in $\Hone$ which satisfy 
\beq \label{e3.5}
\nm{u}_{\pal} \leqs 1 \qquad \mbox{and} \qquad \  \ang{\gamma u, \gamma s_l}_{\bdy} \eqs 0 \quad 
\mbox{for} \quad 0  \leq l \leq k-1.
\eeq

The existence and some properties of such eigenfunctions are described in sections  6 and 7 of \cite{Au1} 
for a more general system. 
In particular, that analysis shows that  each $\delta_j$ is of finite multiplicity and $\delta_j \ra \infty$ as $j \ra \infty$; 
see Theorem 7.2 of   \cite{Au1}.   
The maximizers not only are $\pal -$orthonormal but they also satisfy
\barr 
\IOm   \nabla s_k \cdot  \nabla s_l  \; dx \eqs  \IMby \,  s_k \, s_l  \, \dsg \eqs 0 \qquad 
\mbox{for } \quad k \ne l. \label{e3.7} \\
\IOm  | \nabla s_k |^2  \; dx \eqs \frac{\delk}{1+\delk} \qquad \text{and} \quad 
 \IMby \,  |\gamma s_k|^2   \, \dsg \eqs \frac{1}{1+\delk} \quad \mbox{for} \quad k \ge 0.  \label{e3.8}
\earr
Recently  Daners \cite{D2}  corollary 4.3 has  shown that, when $\Om$ is a Lipschitz domain, then the Steklov eigenfunctions are 
continuous on $\Omb$. 

The analysis in this paper is based on the fact that the harmonic Steklov eigenfunctions on $\Om$ can be chosen 
to be an orthonormal  basis of both  $\Harm$ and of $\Ltby$.  
Note, however,  that these functions are not $L^2-$orthogonal on $\Om$ - so the following constructions are 
quite different to  results from the much studied theory of  Bergman spaces.

Let $\Sc := \{s_j : j \geq 0 \}$ be the maximal family of  $\pal -$orthonormal  eigenfunctions constructed inductively as above. 
Given $u \in \Hone$, consider  the  series
\beq   \label{e3.9}
P_H \, u(x) \deqs \sum_{j=0}^{\infty} \quad [u,s_j]_{\pal} \  s_j(x).
\eeq
This  series will be called a {\it harmonic Steklov expansion} and represents a harmonic function on $\Om$.
For any region $\Om$ the first term in this expansion is the mean value of the function $u$ on the boundary $\bdy$. 
Then theorem 3.1 of \cite{Au3} may be stated as 

\begin{thm} \label{T3.1}
Assume  $\Om, \bdy$ satisfy (A), then $\Sc$ is an orthonormal basis of $\Harm$.
 $P_H$ defined by (\ref{e3.9}), is the $\pal -$orthogonal projection of $\Hone$ onto $\Harm$.  \end{thm}
 \vspace{-1em}
\bpf
This follows from standard results about orthogonal expansions and theorem 
7.3 of \cite{Au1} which says  that $\Sc$ is a  maximal orthonormal subset of $\Harm$. 
\epf

In this paper, it is more convenient to use the Steklov eigenfunctions normalized by their boundary norms.
Define the functions $\sjt(x) \deqs \sqrt{1+\delj} \ \sj(x)$ for $j \geq 0$. 
From \eqref{e3.8}, these satisfy
\beq \label{e3.10}
 \IMby \tilde{s_j} \, \tilde{s}_k \dsg \eqs 0 \xwhen j \neq k \xand  \IMby {\tilde{s_j}}^2  \dsg \eqs 1. 
   \eeq
 These Steklov eigenfunctions are said to be {\it boundary normalized} and the associated set 
 $\Sct \deqs \{\sjt : j \geq 0\} $ is an orthonormal basis of $\Ltby$. See theorem 4.1 of \cite{Au2}.

When $h \in \Ltby$, then  the Riesz-Fischer theorem and the fact that $\Sc$ is an orthonormal  basis of
$\Ltby$ implies that 
 \beq   \label{e3.11}
h(x) \deqs \sum_{j=0}^{\infty} \ c_j \, \sjt(x)   \xwith c_j \deqs \ang{h,\sjt}_{\bdy}  \eeq
in an $L^2-$sense. 
The first term in this series is  the mean value of h on $\bdy; \ \hav \deqs   \ang{h,1}_{\bdy}$.
The analysis in \cite{Au2} shows that the harmonic extension of this function to $\Om$ has the same 
representation for any $x \in \Om$  and has finite energy on $\Om$ if and only if
\beq \label{RF}
\sum_{j=0}^{\infty} \ (1+\delj) \ |c_j|^2 \  < \ \infty.
\eeq

More details of this are to be found in \cite{Au2} where  \eqref{RF} is used as the definition for  $h$ to be in $\Halfby $.  
% That is \eqref{e3.11}  provide a spectral representation of harmonic functions on $\Om$
 This  Steklov representation of the harmonic function $h$ only  involves  integrals  of the boundary values of  $h$.   

Here our interest is in the approximation of harmonic functions by finite sums of this series.  
 When $h \in \Harm$ is an $H^1-$harmonic function on $\Om$, then the M-th Steklov approximation of $h$
  is the function 
  \beq  \label{e3.15}
  h_M(x) \deqs \hav \pls \sum_{j=1}^M \  \ang{h,\sjt}_{\bdy} \ \sjt(x)  . \eeq
 % Here $\hav \deqs   \ang{h,1}_{\bdy}$ is the mean value of h on $\bdy$ and we'll often write
  %$\hjhat $  for the j-th Steklov coefficient  $\ang{h,\sjt}_{\bdy} $ in this expression. 
  
 The essential theorem about these approximations is the following.
 \btm \label{T3.3}
 Assume (A) holds, $h \in \Harm$  and $h_M$ is defined by \eqref{e3.15}. Then $h_M$ converges strongly to $h$ in
 $\Hone$ and $\Harm$ and it converges uniformly to $h$ on compact subsets of $\Om$.
 \etm
 \vspace{-3em}
\bpf
The fact that $h_M$ converges strongly to $h$ in $\Harm$ and $\Hone$ follows from the fact that $\Sc$ is an
orthonormal basis of $\Harm$. 
This sequence also converges  uniformly on compact sets from  theorem X, chapter X of Kellogg  \cite{K}.
The theorem in Kellogg  was only proved for regions in $\Rte$ - but a similar  proof  holds in  all dimensions $N \geq 2$. 
\epf

When $\Om$ is a disc in $\Rtwo$, or a ball in $\RN, \ N \geq 3$ these approximations have been 
extensively studied by  Fourier analysts. 
The Steklov eigenfunctions are $r^k \cos{k\theta}, r^k \sin{k \theta} $ for discs and $r^k Y_{kl}(\theta, \phi)$ on a ball in $\Rte$.
Here $Y_{kl}$ is a spherical harmonic function with k,l   integers.
Its worth noting that the usual mean value theorem for harmonic functions on  a disc or a ball follows from this theorem
as $\sj(0) = 0$ for all $j \geq 1$.  

%The following sections will investigate these approximations in much more detail for the special case where 
%$\Om$ is a rectangle in the plane. 

These eigenproblems have a scaling property.  
Given a region $\Omo \subset \RN$, define $\OmL := \{ Lx : x \in \Omo \}$ with $L > 0$. 
When  $h$ is a harmonic function on $\Omo$, then the function $h_L(y) := h(y/L) $ will be a harmonic function on $\OmL$.
When $h$ is a harmonic Steklov eigenfunction on $\Omo$ associated with an $\xSvl \ \delta$, then $h_L$ will
 be a harmonic Steklov eigenfunction on $\OmL$ with the $\xSvl \ \delta / L $.
Thus it will generally suffice to study problems for a normalized bounded region $\Omo$; results for scalings of this
region then follow using this change of variables. 
 
\vspace{1em} 
% ####### Section 4 ######
% This is line 380 approx
\section{Harmonic Steklov Eigenproblems on a Rectangle}\label{s4}

In this paper we will concentrate on the analysis of harmonic functions on the rectangle $\Gal := (-1,1) \times (-\alpha, \alpha)
\subset \Rtwo$. 
Here $\alpha$  is called the {\it aspect ratio} of the rectangle and $ \xwlog$ is assumed to be in $(0,1]$.
This suffices for analysis on any rectangle upon scaling as described above and possibly a rotation by $\pi/2$.

Henceforth points in $\Gal$ will be denoted by $(x,y)$ and the boundary $\pal \Gal$ consists of 4 line segments
$\Gamma_j$ with $\Gamma_1, \Gamma_3$ being the vertical lines with $x=1, \ -1$ respectively and  
$\Gamma_2, \Gamma_4$ being the horizontal lines with $y=\alpha, \ -\alpha$ respectively. 
The perimeter $|\pal \Gal| = 4(1 + \alpha)$.

First note that if $s$ is a harmonic Steklov eigenfunction on $\Gal$, then the functions obtained by symmetrizing 
$s$ about the x and y axes will either be 0 or will again be a harmonic Steklov eigenfunction, so it suffices
to look for eigenfunctions with specific symmetries about the axes.   

Consider the eigenfunctions obtained by separating variables and assuming the following symmetries.  \\
Class I: \quad $\xSfn$s that are even in x and y; \\
Class II: \quad  $\xSfn$s that are odd in x and y; \\
Class III: \quad  $\xSfn$s that are even in x and odd in y, and \\
Class IV: \quad  $\xSfn$s that are odd in x and even in y.

The first  eigenfunction $s_0(x,y) \equiv 1$ is in class I and the other (unnormalized)  eigenfunctions have the forms
\begin{eqnarray} 
s(x,y) \deqs \cosh{\nu x} \, \cos{\nu y} &\xwhen& \quad \tan{\nualp} \pls \tanh{\nu} \eqs 0, \label{e4.1} \\
s(x,y) \deqs \cos{\nu x} \, \cosh{\nu y} &\xwhen& \quad \tan{\nu} \pls \tanh{\nualp} \eqs 0. \label{e4.2}
\end{eqnarray}

When $\alpha =1$, the first eigenfunction in class II is $s_3(x,y) = xy$. 
Otherwise the (unnormalized)  eigenfunctions and eigenvalues in this class have the form 
\begin{eqnarray} 
s(x,y) \deqs \sinh{\nu x} \, \sin{\nu y} &\xwhen& \quad \tan{\nualp} \eqs  \tanh{\nu},  \label{e4.3}\\
s(x,y) \deqs \sin{\nu x} \, \sinh{\nu y} &\xwhen& \quad \tan{\nu} \eqs  \tanh{ \nualp} . \label{e4.4}
\end{eqnarray}

Similarly eigenfunctions in class III have the forms
\begin{eqnarray} 
s(x,y) \deqs \cosh{\nu x} \, \sin{\nu y} &\xwhen& \quad \tan{\nualp} \mns \coth{\nu} \eqs  0, \label{e4.5} \\
s(x,y) \deqs \cos{\nu x} \, \sinh{\nu y} &\xwhen& \quad \tan{\nu} \pls \coth{\nualp} \eqs 0 \label{e4.6}
\end{eqnarray}

Finally the eigenfunctions in class IV have the forms
\begin{eqnarray} 
s(x,y) \deqs \sinh{\nu x} \, \cos{\nu y} &\xwhen& \quad \tan{\nualp} \pls \coth{\nu} \eqs 0 \label{e4.7} \\
s(x,y) \deqs \sin{\nu x} \, \cosh{\nu y} &\xwhen& \quad \tan{\nu} \mns \coth{\nualp} \eqs 0. \label{e4.8}
\end{eqnarray}

These formulae show that the eigenfunctions  are products of trigonometric and hyperbolic functions and 
the Steklov boundary conditions enforce conditions on the possible parameter $\nu$.
The equations for $\nu$ will be called the {\it determining equations}  and only need to consider positive solutions $\nuj$
in view of the symmetries present.  

For all  $\nu$, one has  $ \tanh{\nualp} \leqs \tanh{\nu} < 1 < \coth{\nualp} \leqs \coth{\nu}$ and each of  
these terms converges to 1 as $\nu$ increases.
Hence the determining equations approximate either $\tan{\nu}  = \pm 1$ or else $\tan{\nualp}  = \pm 1$
when  $\nu,  \, \nualp$ are large enough.
Thus there are two families of values of these $\nu$ when $\alpha < 1$; these families asymptote either
 to odd multiples of $\pi/4$ or to  odd multiples of ${\pi}/{4 \alpha}$. 

The $\xSvl \mbox{s}$ are \\
(i) \quad $\delta =  \nu \tanh{\nu} \quad $ when $\nu$ is a solution of the equation in  \eqref{e4.1} or \eqref{e4.5}. \\
(ii) \quad $\delta =  \nu \tanh{(\nualp)} \quad$ when $\nu$ is a solution of the equation in  \eqref{e4.2} or \eqref{e4.8}. \\
(iii) \quad $\delta =  \nu \coth{\nu} \quad $ when $\nu$ is a solution of the equation in   \eqref{e4.3} or \eqref{e4.7}. \\
(iv) \quad $\delta =  \nu \coth{(\nualp)} \quad $ when $\nu$ is a solution of the equation in  \eqref{e4.4} or \eqref{e4.8}. 

These  eigenvalues all are  positive and grow asymptotically linearly with the values $\nu$.
For given $\alpha \in (0,1]$, the functions $\nu \tanh{(\nualp)}$ and $\nu \coth{(\nualp)}$ are strictly increasing 
functions of $\nu$ on  $[0,\infty)$ that have smooth inverses $g_1, g_2$ respectively.
These inverse functions also are strictly increasing and satisfy  $g_2(\delta) < \delta <  g_1(\delta) $ for all $\delta > 0$.
In the following,  the parameter $\nu$ rather than  the Steklov eigenvalues $\delta$ will be used for  the  analysis 
to simplify many expressions. 
For $\delta = O(1)$ , however, one has that $g_1(\delta)$ and $g_2(\delta)$ both are approximately equal to $\delta$
so they  are essentially the same parameter. 

There is an asymptotic ordering of  the classes of eigenfunctions. 
When $\alpha  =1$ all the eigenvalues will be multiple eigenvalues except for $\delta_0$ and $\delta_3$.
For all $\alpha$ the large eigenvalues  appear to be double, and then quartic,  though in fact  they differ by 
exponentially small amounts. 

The boundary normalized eigenfunctions $\sjt$ used in section \ref{s3} are required to satisfy \eqref{e3.10}. 
As will be seen our interest is primarily in $\xSfn$s of class I, and the boundary integrals  of $s^2$ are found to be 
\begin{eqnarray}   \label{e4.11}
 I(\alpha, \nu) & \eqs& 2 \, \alpha \, {\cosh}^2{\nu}\,  [ 1 + \sinc \, ( 2 \nualp)] \pls {\cos}^2{(\nualp)}\, [2 + {\nu}^{-1} \sinh{(2 \nualp)}]  \\
 J(\alpha, \nu) &\eqs& 2 \, {\cosh}^2{\nualp}\,  [ 1 + \sinc \, ( 2 \nu)] \pls  {\cos}^2{(\nu)}\, [2 \alpha +  {\nu}^{-1} \,  \sinh{(2 \nualp)}]
 \label{e4.12}
\end{eqnarray}
for the two functions in \eqref{e4.1} and \eqref{e4.2} respectively. 
Here sinc$(\theta) \deqs   {\theta}^{-1} \,  \sin{\theta}$ is the cardinal sine function.

Thus the normalized eigenfunction $s_{1j}$ associated with the j-th  solution ${\nu}_j$ of \eqref{e4.1} is
\beq  \label{e4.13}
\st_{1j} (x,y) \deqs \sqrt{\frac{4(1+\alpha)}{I(\alpha, \nuj)}} \ \cosh{(\nuj x)} \ \cos{(\nuj y)}
\eeq
while that associated with the  j-th  solution ${\nu}_j$ of \eqref{e4.2} is
\beq  \label{e4.14}
\st_{2j} (x,y) \deqs \sqrt{\frac{4(1+\alpha)}{J(\alpha, \nuj)}} \ \cos{(\nuj x)} \ \cosh{(\nuj y)}. 
\eeq
These eigenfunctions will be $L^2-$orthogonal on the boundary and $\pal$-orthogonal on $\Gal$ whenever they
arise from distinct $\nuj$.

A proof that these eigenfunctions provide a basis for the set of all  Steklov eigenfunctions for the case of a square is outlined 
in section 3.1 of \cite{GP}. 
The argument provided there generalizes to the case of a rectangle in a straightforward manner.

The general results of \cite{Au1} and theorem \ref{T3.3} say that Steklov approximations of the form \eqref{e3.15}
of the solutions of various harmonic boundary value problems converge strongly to harmonic functions in $H^1(\Gal)$. 
 This was explored numerically by the second author in his recent doctoral thesis \cite{Cho}. 
 Given that $\Harm$ is also a RKH space, this implies also that there is pointwise convergence of these approximations 
 at interior points of $\Gal$.
 
 For most  of the following sections, our interest is in studying approximations of harmonic functions using only harmonic Steklov 
eigenfunctions corresponding to small eigenvalues $\delj$. 
Numerically our finite element computations found that all of the lowest harmonic Steklov eigenfunctions were represented  
by some of eigenfunctions of separated variables form.
 The surprising result that was observed is that this convergence is extremely rapid at the center of the rectangle.
 So the value of a solution at the origin may be estimated quite accurately  by evaluating just  a few boundary integrals.
 The following sections provide an analysis of this phenomenon.   

\vspace{1em}
% ####### Section 5 ######
% This is line 620 approx
\section{Central Value Formulae on a Square.}\label{s5}

In this section the results of the previous two sections will be combined to give explicit approximations  
for the value of a harmonic function $h$ at the origin  in terms of boundary integrals of $h$  on the square $G_1$
with sides of length 2. 
The formulae say that $h(0,0)$ is the mean value of $h$ on the boundary plus some simple integrals
with coefficients that converge geometrically to zero.

First note that the Steklov eigenfunctions of classes  II, III and IV all are zero at the origin and that the determining 
equations for the two different eigenfunctions of class I will be the same when $\alpha = 1$.
Thus in this section $\nuj$ will denote the j-th strictly positive solution of 
\beq \label{e5.1}
\tan{\nu} \pls \tanh{\nu} \eqs 0 \eeq
The set of all solutions of this equation is denoted $\Kc_1$.

%%%%% 
The following table gives the value of the first 6 such zeroes to 8 decimal places, the first differences $\varDelta \nuj:=\nuj-\nu_{j-1}$, 
and the Steklov eigenvalues.  
Henceforth we will only consider harmonic Steklov eigenfunctions of class I so that $\delj$ will be the j-th Steklov eigenvalue
of this class (not over all eigenvalues).  
Note that to six decimal places $\delj = \nuj$ already when $j=4$, and the differences converge rapidly to $\pi$.\\

% ################# Insert table here ##############
\begin{table}[here!]
\begin{center}
\resizebox{16cm}{1.0cm}{
\begin{tabular}{|c||c|c|c|c|c|c|}
\hline
j& 1&2&3&4&5&6 \\ \hline
$\nuj$&  2.36502037&5.49780392& 8.63937983& 11.7809725& 14.9225651& 18.0641578 \\ \hline
$\varDelta \nuj$&    &3.13278355& 3.14157591& 3.14159262& 3.14159265& 3.14159265 \\ \hline
$\delta_{j}$& 2.32363775& 5.49761947& 8.63937929&11.7809724& 14.9225651& 18.0641578\\ \hline
\end{tabular}}
\vspace{0.8em}
\caption{The first 6 points in $\Kc_1$, their first differences and the Steklov eigenvalues.}
\end{center}
\end{table}

For larger j, the solutions are approximated  by 
\beq \label{e5.3}
 \nuj \eqs \left(\frac{3}{4} \pls j \right) \,  \pi \quad \mbox{and} \quad 0.999\, \nuj < \  \delj \ < \ \nuj \eeq
as the second term in  \eqref{e5.1} is essentially 1.

The m-th Steklov approximation of the value of $h(0,0)$ is, from \eqref{e3.15} 
  \beq  \label{e5.5}
  h_m(0,0) \deqs \hav \pls \sum_{j=1}^m \  [\ang{h,\soj + \stj}_{\bGo}] \ \soj(0,0)   . \eeq
  Here we have used the fact that $\soj(0,0) = \stj(0,0)$ to simplify the expression. 
  This involves the mean value of $h$ on the boundary together with m further boundary integrals.
  
  Numerically the values of  of $\soj(0,0) = \stj(0,0)$ for small j were found to be

\begin{table}[here!]
\begin{center}
\resizebox{15cm}{1cm}{
\begin{tabular}{|c||c|c|c|c|c|c|}
\hline
j& 1&2&3&4&5&6 \\ \hline
$\st_{1j}(0,0)$&    0.36925721&    0.016382475  &    $7.079865\times 10^{-4}$   &    $3.0594874\times 10^{-5}$  &    $1.3221244\times 10^{-6}$ &    $5.7134174\times 10^{-8}$\\ \hline
\end{tabular}}
\vspace{0.8em}
\caption{The values of the first 6 Steklov eigenfunctions of class I at (0,0)}
\end{center}
\end{table}

 Let  $\tsoj, \tstj$ denote the unnormalized Steklov eigenfunctions of \eqref{e4.1} and \eqref{e4.2}, then this approximation 
  becomes
    \beq  \label{e5.6}
  h_m(0,0) \deqs \hav \pls \sum_{j=1}^m \  c_j \ \int_{\pal G_1} \, h \, (\tsoj + \tstj) \dsg.
   \eeq
 with $c_j := 1/I(1,\nuj)$.   The actual values of $c_j$  and $c_{j}/c_{j-1}$ found computationally were
 
  %############### Inset table here ###################

\begin{table}[here!]
\begin{center}
\resizebox{15cm}{1cm}{
\begin{tabular}{|c||c|c|c|c|c|c|}
\hline
j& 1&2&3&4&5&6 \\ \hline
$c_j$ & $1.7043861\times 10^{-2}$& $3.35481862\times 10^{-5}$&$ 6.26556108\times 10^{-8}$& $1.17005787\times 10^{-10}$&$ 2.18501606\times 10^{-13}$&
$ 4.08039237\times 10^{-16}$\\ \hline
 $c_{j}/c_{j-1}$&  &$1.9683443 \times 10^{-3}$&$1.8676303  \times 10^{-3}$&$1.8674431 \times 10^{-3}$&
 $1.8674427  \times 10^{-3}$& $1.8674427  \times 10^{-3}$\\ \hline
\end{tabular}}
\vspace{0.8em}
\caption{The values of the first 6 coefficients in \eqref{e5.6}.}
\end{center}
\end{table}

These values appear to  converge rapidly to zero and a logarithmic graph shows that the convergence is geometric. 
This may be  proved using  the following estimates.

\btm \label{T5.1}
When $I(1, \nu)$ is defined by \eqref{e4.11}  then, for $j \geqs 1,$
\beq \label{e5.7}
0 < c_j \,  < 2.56 \,  \exp{(-2 \, \nuj)} \xand 0 < \soj(0,0) \eqs  \stj(0,0) \leqs 4.53 \,  \exp{(-\nuj)}. \eeq
\etm
 \vspace{-2em}
\bpf
The coefficients $c_{j}$ are strictly positive as the $I(1,\nu_{j})$ are.
 To obtain a lower bound, note from (\ref{e4.1})
$|\tan \nu_j|=\tanh \nu_j \leq 1$. Hence $\cos^2 \nu_j \geq \frac{1}{2}$ and the properties of $\sinc$ yield that $1+\sinc 2\nu \geq 0.782$. Substitute in (\ref{e4.11}) then
\begin{align}
%%####################### Please check#############
I(1,\nu_j) &\geq 1.564 \cosh^2 \nu_j + 1 +\frac{\sinh 2\nu_j}{2 \nu_j}\\
&\geq 0.391\ e^{2\nu_j}
\end{align}
Thus $c_j \leq 2.56 \ e^{-2\nu_j}$ and $s_{1j}(0,0) \leq 4.53\ e^{-\nu_j}$ as claimed.
\epf

  For large $j$, this and  \eqref{e5.3} implies that upper bounds for successive coefficients $c_j$ in these approximations decay 
  by a factor of  $e^{-2 \pi} = 1.8674428 \times 10^{-3}$.
  This decay rate is very close to the actual computations for $ j > 2.$ 
  
The extra terms in this expression may be regarded as  {\it correction terms} to the usual mean value theorem that account 
for the fact that   our region is a rectangle rather than a disc. 
Since the $I_{1,\nu}, J_{1,\nu}$ grow exponentially, a computable error estimate may be obtained for this approximation. 
When $h$ is harmonic on $G_1$,  this leads to the following pointwise error bound  for approximating
$h(0,0)$ by $h_m(0,0)$.
\begin{thm} \ ({\bf Harmonic Central Value}) \label{T5.2}
Assume h is a $H^1-$harmonic function on $G_1$, then for $m \geq 1$ there is a $C_m > 0$ such that 
 \beq \label{e5.9}
|h(0,0) \mns h_m(0,0)| \leqs C_m \, e^{-{\nu}_m}  \ \|h||_{\pal G_1}. \eeq
 When $m \geq 3$ then $C_m < 0.41$. \etm
  \vspace{-3em}
\bpf
The Steklov series expansion for $h(0,0)$ converges pointwise at $(0,0)$ from properties of this RKHS and each term satisfies
\begin{equation}
|<h,s_{1j}+s_{2j}>_{\partial G_{1}}s_j(0,0)| \leq 2{\nm{h}}_{\partial G_{1}} s_{1j}(0,0) \leq 9.06 \, e^{-\nu_j} \, {\nm{h}}_{\partial G_{1}}
\end{equation}
upon using the normalization of the eigenfunctions. The boundary norm is the mean $L^2-$norm.
When (\ref{e5.3}) is used,
\begin{align}
 |h(0,0) \mns h_m(0,0)| &\leq 9.06 \, {\nm{h}}_{\partial G_{1}} \, e^{-\nu_{m+1}} \, \sum_{j=0}^{\infty} \, e^{-j\pi}\\
%  &=9.06e^{-\nu_{M+1}} \left( \frac{e^{\pi}}{e^{\pi}-1}\right) {\nm{h}}_{\partial G_{1}} \label{eqth5.2}\\
  &=9.06 \, \left(\frac{e^{-\nu_m}}{e^{\pi}-1} \right) \, {\nm{h}}_{\partial G_{1}} \\
  &\leq 0.41e^{-\nu_m}{\nm{h}}_{\partial G_{1}} \hspace {2cm}\text{as} \ e^{\pi}-1=22.1407
 \end{align}
\epf

That is, these Steklov approximations  converge exponentially to $h(0,0)$  with successive error estimates decreasing by 
factors of $e^{-\pi} = 0.04321392$ when \eqref{e5.3} is used. 
This is consistent with the observed  decrease for the numerical values of the coefficients seen in Table 3.
The following table gives upper bounds on the relative error for some small values of m.
 When m=3, it is an upper bound on $0.41e^{-\nu_3}$

\begin{table}[here!]
\begin{center}

\begin{tabular}{|c||c|c|c|}
\hline
m& 1&2&3 \\ \hline
relative error &    0.039 &    $1.7 \times 10^{-3}$  &    $7.26\times 10^{-5}$  \\ \hline
\end{tabular}
\vspace{0.8em}
\caption{Relative error coefficient for the central value.}
\end{center}
\end{table}

\vspace{1em}
% ####### Section 6  ######
% This is line 700 approx
\section{Central  Value Approximations  on  Rectangles.}\label{s6}

Here  we shall describe the corresponding results when the domain is a rectangle $\Gal$ with $\alpha < 1$. 
First note that attention may be restricted to functions of class I as the other eigenfunctions are zero at the origin. 
In this case the two determining equations differ and the values of the functions $\soj$ and $\stj$
at the origin will differ.

 Let $\nu_{j}^{(1)}$ be the j-th the strictly positive solution of 
 \beq \label{e6.1}
 \tan \alpha \nu = -\tanh \nu
 \eeq
 and   $\nu_{j}^{(2)}$ be that of
 \beq\label{e6.2}
 \tan \nu = - \tanh \alpha \nu
 \eeq
 For moderate and large $\nu$, these are given by
 \beq\label{e6.3}
\nu_{j}^{(1)}=\left(\frac{3}{4}\pls j\right) \ \frac{\pi}{\alpha} \qquad \text{and} \qquad  \nu_{j}^{(2)}\approx \ \alpha \  \nu_{j}^{(1)}\ < \  \nu_{j}^{(1)}
 \eeq
As before the  m-th Steklov approximation to $h(0,0)$ is
\beq\label{e6.4}
h_m(0,0)= \overline{h}\pls \sum_{j=1}^{m} \ \IbG \, h \left[c_{1j}  \tsoj \pls c_{2j} \tstj \right]  \dsg. \eeq
Here  $c_{1j}=1/{I(\alpha,\nuj^{(1)})}$ and  $c_{2j}=1/{J(\alpha,\nuj^{(2)})}$.
The following result provides bounds on the coefficients in the {\it correction terms} needed to accurately estimate the central value of $h$.

\btm \label{T6.1}
When $I(\alpha,\nu)$,   $J(\alpha,\nu)$ are defined by (\ref{e4.11}), (\ref{e4.12}) and $c_{1j}, c_{2j}$ as above then
\begin{align}
  & (i) \quad 0\ <\ c_{1j}\ <\ \frac{2.56}{\alpha}\ e^{-2\nuj^{(1)}} \quad  \text{and also}\quad c_{1j}\  <\  4\ \nuj^{(1)} \  e^{-2 \alpha \nuj^{(1)}}\\
   \text{and} \quad & (ii) \quad 0\ <\ c_{2j}\  < \ 2.56 \ e^{-2\alpha \nuj^{(2)}}
\end{align}
 \vspace{-2em}
 \bpf
 From (\ref{e4.11}) using (\ref{e4.1}) as before, one finds
 \begin{align}
  I(\alpha,\nuj^{(1)})&\geq 1.564\  \alpha \ \cosh^2 \nuj^{(1)} \pls1\pls \frac{\sinh 2\alpha \nuj^{(1)}}{2\nuj^{(1)}}\\
  &\geq 0.391 \  \alpha\    e^{2\nuj^{(1)}}
 \end{align}
 This yields the first part of (i). When $\alpha$ is small this goes to zero so for this case use  
 \beq
  I(\alpha,\nuj^{(1)}) \geq \frac{e^{2\alpha\nuj^{(1)}}}{4\nuj^{(1)}}\left[1\pls 1.564\ \alpha\ \nuj^{(1)}e^{2(1-\alpha)\nuj^{(1)}}\right]
\eeq
which yields  the second inequality in (i). This proves that there is exponential decay even with small $\alpha$. 

 To prove (ii), note that from (\ref{e4.12})
 \begin{align}
J(\alpha,\nuj^{(2)}) &\geq 1.564 \cosh^2 \alpha \nuj^{(2)} \pls 2\alpha\\
&\geq 0.391\  e^{2 \alpha \nuj^{(2)}}
 \end{align}
Hence (ii) follows.
 \epf
\etm
In general one has that $\alpha \nuj^{(2)}\ <\ \nuj^{(2)}\ < \ \alpha \nuj^{(1)} \ < \ \nuj^{(1)}$ so the coefficients $c_{2j}$ 
typically are larger than the $c_{1j}$ but both decrease exponentially as j increases.
Combining \eqref{e6.3} with this theorem one sees that only  good approximations  of $h(0,0)$ may be obtained with 
just a few additional boundary integrals.
The proof of theorem \ref{T5.2} is easily  modified to obtain a similar estimate here and the terms converge geometrically
with factors smaller than or equal to $e^{-\pi / \alpha}$.

\vspace{1em}
% ####### Section 7  ######
% This is line 840 approx
\section{Central Value Formulae for Robin Harmonic Problems.}\label{s7}

The results in the previous sections involved solutions of Dirichlet problems for Laplace's equations. 
That is they provide formulae for $h(0,0)$ when h is harmonic and known on the boundary.

 The analysis can be extended to estimate the central value $h(0,0)$ of a harmonic function $h$ when Robin
  boundary conditions of the form
 \beq \label{e7.1}
 (1-t)D_{\nu}h \pls th \eqs \eta , \quad \text{on}\quad \partial G_{\alpha}
 \eeq
are prescribed. 
Here $t\in (0,1]$ is a constant and it suffices that $\eta = L^{p}(\pal \Om)$ for some $p>1$. 
When this holds then the linear functional
\beq
b(u) \deqs \int_{\pal G_{\alpha}}  \eta \ u \ \dsg
\eeq
is continuous on $H^{1}(G_{\alpha})$ and the solution of Laplace's equation subject to (\ref{e7.1}) is given by
\beq \label{e7.5}
h(x,y) \eqs \frac{\etab}{t} \  \pls \sum_{j=1}^{\infty} \frac{<\eta,\sjt>_{\pal \Om}}{(1-t)\delj + t} \ \sjt(x,y).
\eeq

This series converges in $H^{1}(\Gal)$ with the limit being  a $C^{\infty}$ function on $G_{\alpha}$.
The convergence is uniform on compact subsets of $G_{\alpha}$.
 In particular one sees  that 
 \beq
 h(0,0) \eqs \frac{\etab}{t} \  \pls \sum\limits_{j=1}^{\infty} \ \frac{<\eta,\sjt>_{\pal \Om}}{(1-t)\delj + t} \ \sjt(0,0).
 \eeq
 This summation is over all the Steklov eigenfunctions.
 Since the only eigenfunctions that are non-zero at the origin are  those of class I, the summation may be restricted to this class.
 
 Then the m-th Steklov approximation of $h(0,0)$ is
 \beq\label{e7.2}
 h_m(0,0)=\frac{\etab}{t} \pls \sum_{j=1}^m \ \frac{\sjt(0,0) }{\delj + t(1-\delj)}\  \int_{\pal G_{\alpha}}\eta\ [c_{1j} \,\tsoj  \pls c_{2j} \, \tstj \, ] \ \dsg   
 \eeq
with $c_{1j}, c_{2j}$ as below equation (\ref{e6.4}).
Note that the integrals here are the same as in \eqref{e6.4} - and this reduces to that equation when t=1. 
Thus these coefficients decay exponentially as described in theorem \ref{T6.1} and these partial sums converge rapidly to the
actual value of $h(0,0)$ - which will be  close to $t^{-1} \etab$. 

When $t = 0$ in the boundary condition \eqref{e7.1}, we have  a Neumann harmonic boundary value problem that 
has  solutions if and only if $\ \etab = 0$.
When $t = 0$ and $\etab =0$, the problem has a 1-parameter family of solutions and the solution with mean value 0 on $\Gal$
is given by the infinite sum in \eqref{e7.5}.  
When $\etab = 0$, these formulae show that  as $t\searrow 0^+$ the solutions of  Robin problems converge pointwise to this mean 
zero solution of the Neumann problem

% % %
% ########    Bibliography  #########
%

\end{document}